\documentclass[11pt]{amsart}
\usepackage{amsfonts,amsmath,amssymb}
\usepackage{hyperref}
\begin{document}

\newcommand{\J}{{\rm j}}
\newcommand{\Q}{{\mathbb Q}}
\newcommand{\Z}{{\mathbb Z}}
\newcommand{\F}{{\mathbb F}}
\newcommand{\bF}{{\overline{\mathbb{F}}_p}}
\newcommand{\Ell}{{\rm Ell}}
\newcommand{\Sl}{{\rm{Sl}_2}}
\newcommand{\rs}{{\rm S}}
\newcommand{\M}{{\mathbb M}}
\newcommand{\sK}{{\sf K}}
\newcommand{\ie}{\textit{i}.\textit{e}.\,}
\newcommand{\eg}{\textit{e}.\textit{g}.\,}
\newcommand{\cf}{{\textit{cf}.\,}}
\parindent=0pt
\parskip=6pt

\title{What $\Ell$ sees that $K$ doesn't (when $p > 3$)}

\author[Jack Morava]{Jack Morava}

\address{Department of Mathematics, The Johns Hopkins University,
Baltimore, Maryland 21218}

  
\thanks{Thanks to the organizers and participants at the May 2020 Perimeter 
Institute conference on elliptic cohomology and physics, and especially to
Andrew Baker, Nitu Kitchloo, and Ken Ono for conversations and help with
the content of this note.}  
      
\email{jack@math.jhu.edu}

\date{June 2020}

\begin{abstract}{We use Andrew Baker's analysis of the cofiber of the 
endomorphism
\[
E_{p-1} : S^{2(p-1)} \Ell_p \to \Ell_p
\]
of the $p$-adic elliptic spectrum ($p>3$) to present its completion
away from the locus of ordinary elliptic curves as a sum of roughly 
$p/12$ copies (indexed by supesingular elliptic curves) of $p$-adic 
lifts of the height two mod $p$ cohomology theory $K(2)$. See the
recent paper of Zhu Yifei \cite{8} for a much deeper exploration
of the topics considered in this note.}\end{abstract}

\maketitle \bigskip

{\bf 1} An elliptic curve over a field of characteristic $p$ is said to 
be supersingular if it has no $p$-torsion geometric points; or,
more precisely, if its group (-valued functor) of $p$-torsion points
is connected (and infinitesimal). Such curves have unusually large
groups of automorphisms: their associated formal group laws have height
two. That is, in the formal multiplication-by-$p$ endomorphism 
\[
T \mapsto [p]_E(T) = v_1(E) T^p + v_2(E) T^{p^2} + \cdots \;,
\]
$v_1(E) = 0$, while $v_2(E)$ is a unit. 

Homotopy theorists know that, for primes greater than three and away from 
the locus of supersingular curves in the elliptic moduli stack, 
elliptic cohomology and topological complex $K$-theory are Bousfield 
equivalent - that is, they have the same class of acyclic spaces - and 
therefore see pretty much the same things \cite{1}. The corresponding 
equivariant assertion is not so well understood, but I believe it also 
holds for Borel equivariant versions of these cohomology theories. 

{\bf 2 Recollections} The ring $\Z_p$ of $p$-adic integers can be identified
with the ring of Witt vectors over $\F_p$; more generally, $W(\F_{p^n})$ 
is the valuation ring of the unramifield extension $\Q_{p^n}/\Q_p$ defined
by adjoining a primitive $(p^n-1)$th root of unity.  

Let $M = \Z[E_4,E_6]$ be the ring of integral modular forms for $\Sl(\Z)$; 
following \cite{6}, we have
\[
12^3 \Delta = \eta^{24} = E^3_4 - E^2_6 \; {\rm and} \; j = E^3_4\Delta^{-1} 
\;.\]
Elliptic curves are classified (up to isomorphism, over an algebraically closed
field) by their $j$-invariants. If the curve $E : y^2 = x(x-1)(x-\lambda)$ is
written in Legendre's form, then
\[
j = 2^8 \frac{(\lambda^2 - \lambda +1)^2}{\lambda^2 (\lambda - 1)^2}
\]
and Deuring showed that $E$ is supersingular iff 
\[
\sum_{0 \leq k \leq (p-1)/2} \binom{(p-1)/2}{k}^2 \lambda^k \equiv 0
\]
modulo $p$. A theorem of Deligne identifies $v_1$ with $E_{p-1}$ modulo $p$ 
when $p > 3$, and Gross and Landweber \cite{3} show that 
\[
v_2 \equiv (-1)^{\frac{p-1}{2}} \Delta^{\frac{p^2-1}{12}} \equiv  
(-1)^{\frac{p-1}{2}}12^{-\frac{p^2-1}{4}} \cdot \eta^{2(p^2-1)} \;. 
\]

{\bf 3 Definition} Let 
\[
\M^0_p := W(\F_{p^2})[e_4,e_6]/(3e_6^2 = 1 - e_4^3) \;;
\]
then 
\[
E_4 \mapsto 12\eta^8e_4 ,\; E_6 \mapsto 24\eta^{12}e_6,\; j \mapsto e_4^3 
\]
defines a flat extension 
\[
\Delta^{-1} M \to \M^0_p[\eta^{\pm 1}]
\]
of the ground ring for Landweber-Ravenel-Stong elliptic cohomology, and thus 
a $p$-adic version of elliptic (co)homology at primes greater than 
three\begin{footnote}{Exercise: 3 has a square root mod $p>3$ iff $p \equiv
\pm 1$ mod 12.}\end{footnote}. 

A supersingular elliptic curve $E$ over $\bF$ can be shown to be defined 
over $\F_{p^2}$, implying that $j(E) \in \F_{p^2}$ has Teichm\"uller lift 
$\j(E)$ in $W(\F_{p^2})$\begin{footnote}{Andrew Ogg observes that the only 
primes for which all supersingular $j$-values lie in the prime field are 
exactly those primes occurring in the order of the Monster simple group
\cite{2}.}\end{footnote}. The polynomials
\[
\rs_p(j) = \prod_{E \; {\rm supersingular}/ \bF}(j - \J(E)) \in W(\F_{p^2})[j]
\]
have simple roots, and are of degree 
\[
\sigma(p) := 1 - \epsilon(p) +  \lfloor p/12 \rfloor \;,
\]
where $\epsilon(p) = \pm 1$ if $p \equiv \pm 1$ mod 12, and is otherwise zero.
These polynomials are tabulated for primes less than 53 in \cite{6}(\S 2.8.3). 

{\bf 4} The ring extension 
\[
\M^0_p \to \M^0_{\hat{S}p} := (\M^0_p)_{\widehat{(\rs_p)(j)}}
\]
defined by completion with respect to the ideal generated by $\rs_p(j)$ 
is again flat, defining a version
\[
\M^*_{\hat{S}p} := \M^0_{\hat{S}p} \otimes_M {\rm Ell}^* 
\]
of classical $p$-adic elliptic cohomology completed away from the locus of 
ordinary elliptic curves (\ie whose formal group laws have multiplicative
reduction, defining cohomology theories which, like Tate $K$-theory 
\cite{4}(Th 5.1), are equivalent to classical $K$-theory, but with 
nonstandard complex orientations).

As suggested in \cite{5}{\S 2.2.1}, see further \cite{7},  
\[
\M^*_{\hat{S}p} \cong \oplus \sK^*(\Q_{p^2})
\]
splits as a sum of $\sigma(p)$ lifts to $W(\F_{p^2})$ of the 
$\F_p$-module-valued cohomology functors $K(2)$. \bigskip

\bibliographystyle{amsplain}

\end{document}